\documentclass[12pt,reqno, oneside]{amsart}
\usepackage{pdfsync}  
\usepackage[height=9 in]{geometry}
\usepackage{hyperref}
\begin{document}
\title{On a polynomial congruence for Eulerian polynomials}
\author{Ira M. Gessel$^*$}
\address{Department of Mathematics\\
   Brandeis University\\
   Waltham, MA 02453}
\email{gessel@brandeis.edu}
\date{January 18, 2021}
\thanks{Supported by a grant from the Simons Foundation (\#427060, Ira Gessel).}
%\begin{abstract}
%Abstract goes here.
%\end{abstract}

\maketitle
\thispagestyle{empty}

Yoshinaga \cite[Proposition 5.5]{yoshinaga} proved, using arrangements of hyperplanes, the polynomial congruence for Eulerian polynomials
\begin{equation}
\label{e-cong}
A_n(t^m) \equiv \left(\frac{1+t+\cdots+t^{m-1}}{m}\right)^{n+1}A_n(t) \mod (t-1)^{n+1}.
\end{equation}
Here the Eulerian polynomials $A_n(t)$ may be defined by the generating function
\begin{equation}
\label{e-gf}
\sum_{n=0}^\infty \frac{A_n(t)}{(1-t)^{n+1}}\frac{x^n}{n!}=\frac{1}{1-te^x}.
\end{equation}
A simpler proof, using roots of unity, was given by Iijima et al.~\cite{iijima}.
We give here a very simple proof based on the generating function \eqref{e-gf}.

Since $1+t+\cdots+t^{m-1}=(1-t^m)/(1-t)$, the congruence \eqref{e-cong} is equivalent to the statement that the denominator of 
\begin{equation}
\label{e-dif}
m^{n+1}\frac{A_n(t^m)}{(1-t^m)^{n+1}} - \frac{A_n(t)}{(1-t)^{n+1}}
\end{equation}
is not divisible by $t-1$. 

By \eqref{e-gf} the rational function \eqref{e-dif} is the coefficient of $x^n/n!$ in 
\begin{equation*}
\frac{m}{1-t^m e^{mx}} -\frac{1}{1-te^x}
 =\frac{m}{1-t^m e^{mx}} - \frac{1+te^x +t^2 e^{2x}+\cdots + t^{m-1}e^{(m-1)x}}{1-t^m e^{mx}}.
\end{equation*}
Thus it suffices to show that the denominator of the coefficient of $x^n/n!$ in 
\begin{equation}
\label{e-mj}
\frac{1-t^j e^{jx}}{1-t^m e^{mx}}=\frac{1+te^x +\cdots + t^{j-1}e^{(j-1)x}}{1+te^x +\cdots + t^{m-1}e^{(m-1)x}}
\end{equation}
is not divisible by $t-1$.
We have $1+te^x +\cdots + t^{m-1}e^{(m-1)x} = 1+t+\cdots + t^{m-1} +x P(t,x)$
where $P(t,x)$ is a power series in $x$ with coefficients that are polynomials in $t$. 
It follows that the denominator of the coefficient of $x^n/n!$ in in \eqref{e-mj} is a constant times a power of $1+t+\cdots+t^{m-1}$ and is thus not divisible by $t-1$.

%for BibTex
\bibliography{eulerian}{}
\bibliographystyle{amsplain}

\end{document}